\documentclass{article}
\usepackage{amsmath}
\usepackage{amssymb}
\newtheorem{thm}{Theorem}[section]
\newtheorem{lem}[thm]{Lemma}
\newtheorem{prop}[thm]{Proposition}

\title{Arf characters of an algebroid curve}
\author{V. Barucci\thanks{{\em email} barucci@mat.uniroma1.it}
\and M. D'Anna\thanks{
{\em email} mdanna@dmi.unict.it,
partially supported by \lq\lq Progetto Giovani Ricercatori" of the
University of Catania.}\and R. Fr\"oberg\thanks
{{\em email} ralff@matematik.su.se}} 

\begin{document}
\date{}
\maketitle
\begin{abstract}  

Two algebroid branches are said to be equivalent if they have the same
multiplicity sequence. It is known that  two algebroid branches $R$ and $T$
 are equivalent if
and only if their Arf closures,
$R'$ and $T'$ have the same value semigroup, which is an Arf
numerical semigroup and can be expressed in terms of
a finite set of information, a set of characters of the
branch. 

We extend the above equivalence to algebroid curves
with $d>1$ branches. An equivalence class is described, in this more
general context, by an Arf semigroup, that is not a numerical semigroup,
but is a subsemigroup of $\mathbb N^d$. We express this semigroup in terms
of a finite set of information, a set of characters of the
curve, and apply this result to determine other curves equivalent to a
given one. 
\medskip
 
\end{abstract}

\section{Introduction}
By an algebroid branch we mean a one-dimensional domain of the form
$R=k[[x_1,\ldots,x_n]]/P$, where $x_1,\ldots,x_n$ are indeterminates over
the field $k$ (that we assume to be an algebraically closed field of characteristic zero) and $P$ is a prime ideal in $k[[x_1,\ldots,x_n]]$. 
The integral closure
$\bar R$ of
$R$ is isomorphic to the DVR $k[[t]]$. Thus every nonzero element in $R$
has a value, considering $R$ as a subring of $\bar R$. The set of values of
nonzero elements in $R$ constitute a numerical semigroup $v(R)=S$, i.e. an
additive submonoid of $\mathbb N$, with finite complement in $\mathbb N$,
 and the multiplicity of the ring $R$, $e(R)$, is given by the smallest
positive value in $S$. The blowup of $R=R_0$ is $R_1= \bigcup_{n \geq 0}
(m^n:m^n) = R[x^{-1}m]$, where $m$ is the maximal ideal of $R$ and $x$ is
an element  of smallest value in $m$. The blowup $R_1$ is a local
overring of $R$ and, if $R_{i+1}$ denotes the blowup of $R_i$, then
$R_j=k[[t]]$, for $j>>0$. The multiplicity sequence of $R$ is defined to
be $e_0=e(R_0), e_1=e(R_1), \ldots$.

We define two algebroid branches as  equivalent if they
have the same multiplicity sequence. This equivalence extends the Zariski
equivalence between plane branches (cf. \cite {zar})  to branches
of any embedding dimension and has been studied by several authors
(cf. e.g.  \cite [Definition 1.5.11]{C}). 

As has been proved by the work of Arf \cite{arf}, Du Val \cite {D1},
Lipman \cite {L}, two algebroid branches $R$ and $T$ are equivalent if
and only if their Arf closures,
$R'$ and $T'$ have the same value semigroup, that is an Arf
numerical semigroup. For the convenience of the reader these results are
set out in detail in Section 2. In Section 3 an explanation is given on
how to express an Arf numerical semigroup (which always represents the
equivalence class of a branch $R$) in terms of a minimal finite set of
information, the characters of $R$. The content of Proposition \ref {Arf}
is more or less implicit but not explicitely proved in \cite {D1}. 

The main results of the paper are contained in Sections 4 and 5. In Section
4 we extend the above equivalence to algebroid curves with $d>1$ branches.
As in the one branch case, an equivalence class is described, in this more
general context, by an Arf semigroup, that is not a numerical semigroup,
but is a subsemigroup of $\mathbb N^d$. On the other hand any such
semigroup describes an equivalence class (cf. Theorem \ref {bij}).
Although  the value semigroup $v(R)$ of an algebroid curve $R$ with $d>1$
branches is not finitely generated, the equivalence class of $R$,
represented by an Arf semigroup contained in $\mathbb N^d$ or by a
multiplicity tree with $d$ branches as well, can be expressed in terms of
a finite (but not unique) set of information (cf. Theorem \ref {V}).
Such a finite set, a set of characters of the curve $R$, is useful to
determine subrings of $R$ (for example projections of the curve $R$ on
spaces of lower dimension) that are equivalent to $R$ (cf. Theorem
\ref{subrings}). This is illustrated in the last two examples.

The same invariant that we consider for algebroid curves is studied 
 by A. Campillo and J. Castellanos in \cite {cc} for schemes of
arbitrary dimension, with the name of Valuative Arf Characteristic, as
pointed out by the referee. Anyway the results of this paper are coherent
with \cite {cc}, but not contained in it.

Finally we want to thank the referee for his/her carefull reading and
useful suggestions.
 
\section{Equivalence between algebroid branches}
We keep the hypotheses and notation of the Introduction: let $R$ be an
algebroid branch with value semigroup $v(R)=S$, let $R=R_0, R_1, \ldots$
be the sequence of successive blowups of $R$ and let $e_0=e(R_0),
e_1=e(R_1),
\ldots$ be its multiplicity sequence.
Two algebroid branches are said to be {\it equivalent} if they have the
same multiplicity sequence.
If $\alpha$ is a nonnegative integer, let $R(\alpha)=\{r \in R; v(r)
\geq \alpha\}$. Then $R(\alpha)$ is an ideal of $R$ and an ideal $I$ of $R$
is integrally closed if and only if $I=R(\alpha)$, for some $\alpha$. The
ring $R$ is said to be an {\it Arf ring} if $x^{-1}R(\alpha)$, where
$v(x)=\alpha$, is a ring, for each $\alpha \in v(R)$.
For each ring $R$ there is a smallest Arf overring $R'$, called the {\it
Arf closure} of $R$. Since the Arf closure maintains the
multiplicity of the ring and commutes with blowingup (cf.
\cite [Theorem 3.5]{L}), $(R_i)'=(R')_i$, the multiplicity sequence
of $R$ is the same as the multiplicity sequence of $R'$, i.e. $R$ is
equivalent to its Arf closure $R'$.

It is wellknown that the value semigroup of the Arf closure of $R$,
$v(R')$ is an {\it Arf semigroup}, i.e. a numerical semigroup $S$ such
that $S(s)-s$ is  a semigroup, for each $s \in S$, where $S(s)=\{n \in S;
n \geq s\}$ (cf. \cite [page 8]{Ba-Do-Fo}).

Given an Arf numerical semigroup $S=\{s_0=0<s_1<s_2<\cdots\}$, the
{\it multiplicity sequence of} $S$ is the sequence of differences $s_
{i+1}-s_i$, for $i \geq 0$. The first element of this sequence
$s_1-s_0=s_1$ is called the {\it multiplicity of} $S$.
If $S=v(R')$, the multiplicity sequence of the
ring $R$ is exactly the same as the multiplicity sequence of the
Arf semigroup $S$ (cf. \cite{arf} or \cite[Proposition 5.10]{BDF}).

We recall now how the multiplicity sequence of an Arf semigroup
characte-rizes the semigroup completely.

It follows immediately from the definitions that the multiplicity sequence
$e_0, e_1, \ldots$ of an Arf semigroup is such that, for each $i \geq 0$,
$e_i = \sum_{h=1}^{k}e_{i+h}$, for some $k \geq 1$.
Conversely any sequence $e_0, e_1, \ldots$ of natural numbers such that
$e_n=1$, for $n >> 0$ and, for all $i$, $e_i = \sum_{h=1}^{k}e_{i+h}$, 
for some $k \geq
1$, is the multiplicity sequence of the Arf semigroup $S$ given by the
sums, $S=\{0,e_0,e_0+e_1,\dots\}$. We call such a sequence of natural
numbers a {\it multiplicity sequence}. For example, given the sequence
$6,3,3,3,1,1,\ldots$, the Arf semigroup with such multiplicity sequence
is $S=\{0,6,9,12,15,\rightarrow\}$ (with
$\rightarrow$ we mean that all consecutive integers are contained in the
set). Thus to give an Arf numerical semigroup is equivalent to give its
multiplicity sequence.

Moreover, if we consider $S_1=(S\setminus\{0\})-e_0= 
\{0,e_1,e_1+e_2,\dots\}$ and, for $i \geq 1$, $S_{i+1}=(S_i\setminus 
\{0\})-e_i=\{0,e_{i+1},e_{i+1}+e_{i+2},\dots\}$, we get a chain of Arf
semigroups $S \subset S_1 \subset S_2 \subset \cdots\subset S_n =\mathbb N$ 
such that 
the multiplicity of $S_i$ is $e_i$. 

Recall also that, for each numerical semigroup $S$, there is a smallest
Arf semigroup $S'$ containing $S$, called the {\it Arf closure} of $S$. It
is the intersection of all the Arf semigroups containing $S$. If $s_1,
\ldots,s_n$ are natural numbers, we denote by Arf$(s_1,\ldots,s_n)$ the
smallest Arf semigroup that contains $s_1,\ldots,s_n$. Given an Arf
semigroup $S$, there is a uniquely determined smallest semigroup $N$ such
that the Arf closure of $N$ is $S$. The minimal system of generators
$\{n_1,\cdots,n_h\}$ for such $N$ is called in \cite {RGGB} the {\it Arf
system of generators} for $S$. We also will use this terminology.    
The generators $\{n_1,\cdots,n_h\}$  are called in \cite {arf} and \cite
{D1} the {\it characters} of the ring $R$ (if $v(R')=S$); so we refer to
them as the {\em Arf characters} of $R$ and in the next  section we give an
algorithm to find them.

\medskip

In conclusion two algebroid branches $R$ and $T$ are equivalent if and only if their Arf closures $R'$ and $T'$ have the same value semigroup,
that is an Arf semigroup.
On the other hand each Arf semigroup $S$ is the value semigroup 
of the algebroid branch $R=R'=k[[S]]$. It follows that the equivalence classes
of algebroid branches are in one to one correspondence with Arf numerical semigroups.

In the particular case of two algebroid plane branches $R$ and $T$
it is well known that $R$ and $T$ are equivalent if and only if 
$v(R)=v(T)$(cf. \cite{zar86}, \cite {wa} or \cite[Theorem 2.5]{piane}).
So in the plane case $v(R)=v(T)$ if and only if $v(R')=v(T')$.
Notice however that  $v(R')$ is an Arf semigroup that contains ($v(R)$ and so)  the
Arf closure of $v(R)$, but the inclusion $v(R)' \subseteq v(R')$ may be
strict. A simple example is the following.

\medskip
\noindent{\bf Example} If $R=k[[t^4,t^6+t^7]]$, then
$v(R)=\langle4,6,13\rangle=\{0,4,6,8,10,12,13,14,\\ 16,\rightarrow\}$. Thus
$v(R)'=\langle4,6,13,15\rangle=\{0,4,6,8,10,12,\rightarrow\}$. However
$R'=k+kt^4+k(t^6+t^7)+t^8k[[t]]$, thus
$v(R')=\langle4,6,9,11\rangle=\{0,4,6,8,\rightarrow\}$.

\medskip

Notice also that, given an equivalence class of algebroid branches, or, 
equivalently, given an Arf numerical semigroup, not always there is a plane
 branch in that class: there is a characterization of Arf semigroups
representing a class that contains a plane branch (cf. e.g. \cite[Theorem
3.2]{piane}).

\section{Arf characters of an algebroid branch}

As shown in the previous section the 
equivalence class of an algebroid branch $R$ is 
determined by the value semigroup of its Arf closure $v(R')=S$
and hence is given in terms of finite information:
the set of generators of $S$ or, better, the Arf system of generators of $S$.
The advantage with the Arf system of generators is that we can generalize it to curves with several branches 
and it is a finite set also in that case (as we will see in the last section),
while the semigroup is no more finitely generated.

\medskip

We will now give an algorithm to find the Arf system of generators of an
Arf numerical semigroup.

Let $S$ be an Arf numerical semigroup and let
$e_0, e_1, \ldots$ be its multiplicity sequence. 
The {\it restriction
number} $r(e_j)$ of $e_j$ is defined to be the number
of sums $e_i= \sum_{h=1}^{k} e_{i+h}$ where
$e_j$ appears as a summand. With this terminology we 
have the following result:

\begin{prop}\label{Arf} Let $S$ be an Arf numerical semigroup and let
$e_0, e_1, \ldots$ be its multiplicity sequence. Let $r(e_j)$ 
be the restriction
number of $e_j$. Then the Arf system of
generators for
$S$ is given by the elements
$e_0+e_1+\cdots+e_j$, where $r(e_j)<r(e_{j+1})$.
\end{prop}

To prove this result we need a lemma:
\begin{lem} Suppose $e_0,e_1,\cdots$ is a multiplicity sequence and that
$\{0,e_0,e_0+e_1,\ldots\}$ is its associated Arf semigroup. Then
$r(e_j)<r(e_{j+1})$
if and only if there is no $k<j$ such that $e_k=e_{k+1}+...+e_j$. In this
case $r(e_j)=r(e_{j+1})-1$.
\end{lem} 

\noindent {\bf Proof.} There is only one sum in which $e_{j+1}$ but not
$e_j
$ is a summand, namely $e_j=e_{j+1}+\cdots$. All other sums that contain
$e_{j+1}$ as a summand also contain $e_j$ as a summand. Thus
$r(e_j)\ge r(e_{j+1})$ if and only if there is a $k<j$ such that
$e_k=e_{k+1}+\cdots+e_j$. Otherwise $r(e_j)=r(e_{j+1})-1$.

\medskip
\noindent {\bf Proof of Proposition 3.1.} We argue by induction on the
number of elements
$\neq 1$ in the multiplicity sequence  of $S$. If this number is zero,
i.e. the multiplicity sequence is $1,1,\dots$, then $S= \mathbb N= {\rm Arf}(1)$
and in fact only the first $1$ of the multiplicity sequence has restriction
number (equal to zero) strictly smaller than the next. 

Let
$S=\{0,e_0,e_0+e_1,e_0+e_1+e_2,\dots\}$ be an Arf semigroup. Then
$S_1=\{0,e_1,e_1+e_2,\dots\}$ is also an Arf semigroup and the number of
elements
$\neq 1$ in the multiplicity sequence  of $S_1$ 
is one less than the number of elements
$\neq 1$ in the multiplicity sequence  of $S$. 
By the inductive
hypothesis the minimal Arf system of generators for $S_1$ is given by
$\{h'_j=e_1+\dots+e_j; r'(e_j)<r'(e_{j+1})\}$
(where $r'(e_j)$ denotes the restriction number of $e_j$ as
an element of the multiplicity sequence of $S_1$).
We will denote by $J$ the set of the indices of the elements $h'_j$.

We now compare the restriction numbers $r(e_j)$ of the multiplicity
sequence $e_0,e_1,\dots$ of $S$ with the restriction numbers
$r'(e_j)$ of the multiplicity
sequence $e_1,e_2,\dots$ of $S_1$.
We have $r(e_0)=0$ and, if $e_0=e_1+\dots+e_i$,
$r(e_1)=r'(e_1)+1,\dots,r(e_i)=r'(e_i)+1, r(e_j)=r'(e_j)$,
for $j>i$. 
So, excluding the element $e_0$ (of course $r(e_0)=0<r(e_1)$), we have
$r(e_j)<r(e_{j+1})$ if and only if $r'(e_j)< r'(e_{j+1})$,
for every $j \neq i$.
If $j=i$, then $r(e_i)=r'(e_i)+1$ and $r(e_{i+1})=r'(e_{i+1})$.
Hence, even if  $r'(e_i)=r'(e_{i+1})-1$,
$r(e_i) \geq r(e_{i+1})$. 
It follows that 
$\{h_j=e_0+\dots+e_j;\ r(e_j)<r(e_{j+1})\}= 
\{e_0\}\cup \{e_0+h'_j; j \in J,\ j \neq i\}$.
Hence we need to prove that $S={\rm Arf}(e_0,e_0+h'_j;\ j \in J;\ j\neq i)$
and that for every $j\neq i$, $e_0+h'_j$ is necessary as Arf generator
for $S$.

Since $e_0 =e_1+\dots+e_i\in S_1$,  $S_1={\rm Arf}(h'_j=e_1+\dots+e_j
;\ j \in J)={\rm Arf}(e_0, h'_j=e_1+\dots+e_j,\  j \in J)$.
By \cite [Proposition 16]{RGGB}, $S=\{0\} \cup e_0+S_1={\rm Arf}(e_0, 2e_0,
e_0+ h'_j; \ j \in J)$.

Since $2e_0=e_0+e_1+\dots+e_i$ is clearly superfluous, we get
$S={\rm Arf}(e_0, e_0+ h'_j; \ j \in J;\ j \neq i).$

Finally we notice that, for every $j \in J$, $j \neq i$, the elements 
$e_0+h'_j$ are necessary as Arf generators for $S$.
Otherwise we would have, for $H \subsetneq J\setminus\{i\}$,
that $S={\rm Arf}(e_0, e_0+h'_j; j \in H)$ and, again by 
\cite [Proposition 16]{RGGB}, $S_1={\rm Arf}(e_0, h'_j;\ j \in H)$,
a contradiction against the minimality of the set $\{h'_j;\ j\in J\}$
as Arf system of generators for $S_1$.
 
\medskip
\noindent{\bf Example} If $S=\langle6,9,16,17\rangle=
\{0,6,9,12,15,16,\rightarrow\}$, its
multiplicity sequence is $6,3,3,3,1,1,\ldots$, then $r(e_0)=0$, $r(e_1)=1$,
$r(e_2)=2$, $r(e_3)=1$, $r(e_4)=1$, $r(e_5)=2$, $r(e_6)=2$, $r(e_j)=1$,
for $j \geq 7$. The Arf system of generators for $S$ is $\{e_0=6,
e_0+e_1=9, e_0+e_1+e_2+e_3+e_4=16\}$.

\medskip
\noindent{\bf Remark}
Given an Arf semigroup $S$, 
the cardinality of its minimal Arf system of generators
is said to be the Arf rank of $S$, denoted by Arfrank$(S)$
(cf.\cite{RGGB}).

Now let as above $S=\{0,e_0,e_0+e_1,\dots\} 
\subset S_1=\{0,e_1,e_1+e_2,\dots\} \subset S_2= \{0,e_2,e_2+e_3,\dots\}
\subset \cdots \subset S_n=\mathbb N$ be the chain of Arf semigroups
obtained by $S$.

In the proof of \ref{Arf} it is shown that 
${\rm Arfrank}(S_i)={\rm Arfrank}(S_{i+1})+1$,
if $e_i$ is not in the Arf system of generators of $S_{i+1}$,
and ${\rm Arfrank}(S_i)={\rm Arfrank}(S_{i+1})$,
if $e_i$ is in the Arf system of generators of $S_{i+1}$.
It follows that ${\rm Arfrank}(S) \leq n+1$.

It is not difficult to give examples of Arf semigroups $S$ 
such that the equality holds. It is in fact enough to choose for each $i$ 
the multiplicity $e_i$ of $S_i$ as an element of $S_{i+1}$,
that is not in the Arf system of generators of $S_{i+1}$.

\medskip
\noindent{\bf Example}
We give the construction of an Arf semigroup $S$ with 
$n=3$ and ${\rm Arfrank}(S)=4$.

Let $S_3=\mathbb N={\rm Arf}(1)$ (hence its Arfrank is $1$).
Choose a multiplicity for $S_2$ in $S_3 \setminus \{1\}$,
for example $e_2=3$.
Then $S_2=\{0,3,4,\rightarrow\}={\rm Arf}(3,3+1)={\rm Arf}(3,4)$ and 
${\rm Arfrank}(S_2)=2$.

Choose a multiplicity  for $S_1$ in $S_2\setminus\{3,4\}$,
for example $e_1=5$. Then  
$S_1=\{0,5,8,9,\rightarrow\}={\rm Arf}(5,5+3,5+4)={\rm Arf}(5,8,9)$ and 
${\rm Arfrank}(S_1)=3$.

Finally choose a multiplicity for $S$ in $S_1\setminus\{5,8,9\}$,
for example $e_0=10$. 
Then 
$S_1=\{0,10,15,18,19,\rightarrow
\}={\rm Arf}(10,10+5,10+8,10+9)={\rm Arf}(10,15,18,19)$ and ${\rm Arfrank}(S)=4$.

The semigroup $S$ has multiplicity sequence $10,5,3,1,1,\dots$ and the 
number of elements in
this sequence different to $1$ (or, equivalently, 
the length of the chain $S \subset S_1 \subset \cdots \subset S_n=\mathbb N$)
is $n=3$.

\medskip
\noindent{\bf Example} A class of Arf semigroups $S(k)$ such that 
${\rm Arfrank}(S)=k+1$ and $n=k$ is given by 
$S(k)=\langle 2^k,2^k+2^{k-1},\dots,2^k+2^{k-1}+\dots+2+1 \rangle$.


\section{Equivalence between algebroid curves}
Now we consider {\it algebroid curves}, i.e. one-dimensional reduced rings
$R=k[[x_1,\ldots,x_n]]/I$, where $I=P_1 \cap\cdots\cap P_d$, and where
$P_1,\ldots, P_d$ are prime ideals of $k[[x_1,\ldots,x_n]]$. The rings
$k[[x_1,\ldots,x_n]]/P_j=R/\bar P_j$ are called the branches of the curve.
The integral closure 
$\bar R$ of
$R$ is isomorphic to  $k[[t_1]] \times\cdots \times k[[t_d]]$

In this
case the blowup $R_1$ of $R=R_0$, i.e. $R_1= \bigcup_{n \geq 0}
(m^n:m^n)$, where $m$ is the maximal ideal of $R$, is a semilocal ring.
The sequence of overrings $R=R_0 \subseteq R_1 \subseteq \cdots$ is
defined blowing up at each step the Jacobson radical of the previous ring
and $R_j=\bar R=k[[t_1]] \times\cdots \times k[[t_d]]$, for $j>>0$.

Given a maximal ideal $n_j= k[[t_1]] \times \cdots \times t_jk[[t_j]]
\times
\cdots \times k[[t_d]]$ of $\bar R$ the  {\it branch sequence} of $R$ along
$n_j$ is the sequence of local rings $(R_i)_{n_j \cap R_i}$.

We defined in \cite{BDF} the {\it blowing up tree} of $R$
(cf. also \cite{N}): the nodes at
level
$i$ are the local rings $(R_i)_{n_j \cap R_i}$, $1 \leq j \leq d$. For
$i=0$, $(R_0)_{n_j \cap R_0}=R_0=R$, for each $j$. For $i\geq 1$, $R_i$
has a certain number $r_i$ of maximal ideals, $1 \leq r_i \leq d$, and
$r_i \leq r_h$, if $i \leq h$. A node at level $i$ is connected to a node
at level $i+1$ if and only if the corresponding local rings are in the
same branch sequence along some maximal ideal $n_j$ of $\bar R$, i.e. if
and only if their maximal ideals are one over the other. Recall also that,
for each
$i
\geq 0$,
$R_i$ is the product of its localizations at the maximal ideals, i.e. is
the product of the nodes at level $i$ in the blowing up tree (cf. \cite
[Proposition 3.1]{BDF}).

Since $R \subseteq  \bar R=k[[t_1]] \times\cdots \times k[[t_d]]$, every
nonzero divisor in $R$ has a value in $\mathbb N^d$. The set of values of
nonzero divisors in $R$ constitute a subsemigroup of $\mathbb N^d$. This
semigroup $v(R)=S$ satisfies the following conditions, cf. \cite {BDF}:

(1) If ${\boldsymbol\alpha}=(\alpha_1, \dots , \alpha_d)$ and
${\boldsymbol\beta}=(\beta_1,\dots,\beta_d)$ are elements of $S$, then
$\min({\boldsymbol\alpha},{\boldsymbol\beta})=(\min\{\alpha_1,\beta_1\},
\ldots,\min\{\alpha_d,\beta_d\})\in S$.\\
\indent (2) If ${\boldsymbol\alpha},{\boldsymbol\beta}\in S$,
${\boldsymbol\alpha}\ne {\boldsymbol\beta}$ and $\alpha_i=\beta_i$ for
some $i\in\{ 1,\ldots,d\}$, then there exists ${\boldsymbol\epsilon}\in S$
such that $\epsilon_i>
\alpha_i=\beta_i$ and $\epsilon_j\ge\min\{\alpha_j,\beta_j\}$ for each
$j\ne i$ (and if $\alpha_j\ne\beta_j$ the equality holds).\\
\indent (3) There exists ${\boldsymbol\delta}\in{\mathbb N}^d$ such that
$S\supseteq{\boldsymbol\delta}+{\mathbb N}^d$.

Any subsemigroup of $\mathbb N^d$, satifying the three conditions above is
called a {\it good semigroup}. A good semigroup is {\it local} if ${\bf
0}$ is the only element of the semigroup which has some coordinate equal
to $0$. In fact the semigroup $v(R)$ is local if and only if the ring $R$
is local.
 If ${\boldsymbol\alpha}=(\alpha_1,\ldots,\alpha_d)$ is the minimal positive
value in a local semigroup
$S=v(R)$, then the multiplicity of the algebroid curve $R$ is  given by
$ \alpha_1+\cdots+\alpha_d$ (cf. \cite [Theorem 1]{N}). We define the {\it
fine multiplicity} of
$R$ to be the vector ${\boldsymbol\alpha}=(\alpha_1,\ldots,\alpha_d)$. Since
$v(R/\bar P_j)$, the value semigroup of the $j$-th branch of the curve, is
the projection of $v(R)$ on the $j$-th component of $\mathbb N^d$, it
turns out that
$\alpha_j$ is the multiplicity of the branch $R/\bar P_j$.

Let $U$ be the local ring corresponding to a node of the blowing up tree
of the algebroid curve $R$. Suppose that $U$ occurs in the
$j_1,\ldots,j_r$-th branch sequences of $R$, so that $\bar
U=k[[t_{j_1}]]\times \cdots \times [[t_{j_r}]]$.
(Notice that $U$ is the localization of $R_i$ at its $j$-th 
maximal ideal $m_{i,j}=n_{j_1} \cap R_i=\dots =n_{j_r} \cap R_i$;
hence $j_1,\dots j_r$ are consecutive indices.)
Then $U$ has $r$ minimal
primes $q_{j_1},\ldots, q_{j_r}$. We denote in  this case the fine
multiplicity ${\bf e}(U)$ of $U$ by a vector with $d$ (not $r$)
 components, setting ${\bf e}(U)=(e_1(U),\ldots,e_d(U))$, where
  $e_{j_h}(U)=e(U/q_{j_h})$ and $e_j(U)=0$,  if $j \notin \{j_1,\ldots,
j_r\}$. If we replace the local rings in the tree with their fine
multiplicities, we get the {\it multiplicity tree} of the algebroid curve
$R$. If $j=j_h$, for some $h$, and if the $j_1$-th, \dots, $j_r$-th branches
are glued in a node at level $i$, we denote this node on
the $j$-th branch of the multiplicity tree by
$$
{\bf e}_{(i)}^j = (0,\dots ,0, e_i^{j_1},\dots,e_i^{j_r},0,\dots,0).
$$
 
\medskip
\noindent{\bf Example} The following is a very simple example of an
algebroid curve with two branches.

Let $R$ be the subring $k[[(t^2,u^2),(0,u^3),(t^3,0)]]$ of
$k[[t]]\times k[[u]]$. This is the ring:
$$k[[x,y,z]]/(x^3-z^2,y)\cap(x^3-y^2,z)$$
Blowing up the maximal ideal of $R$ we get
$R_1=k[[(t^2,u^2),(0,u),(t,0)]]$, i.e. the ring
$$k[[x,y,z]]/(x-z^2,y)\cap(x-y^2,z)$$ that is still local. At next
blowup the two branches split and we get the semilocal ring
$R_2=k[[t]]\times k[[u]]$.
Thus the multiplicity tree of $R$ is the following

 \setlength{\unitlength}{0.49mm}
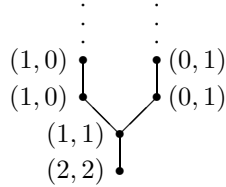
\begin{center}
\begin{picture}(110,55)(26,-75)

\put(55,-72){$(2,2)$}
\put(55,-62){$(1,1)$}
\multiput(45,-52)(0,10){2}{$(1,0)$}
\multiput(88,-52)(0,10){2}{$(0,1)$}
\multiput(75,-70)(0,10){2}{\circle*{2}}
\put(75,-70){\line(0,1){10}}
\put(75,-60){\line(-1,1){10}}
\put(75,-60){\line(1,1){10}}
\multiput(65,-50)(0,10){2}{\circle*{2}}
\put(65,-50){\line(0,1){10}}
\put(85,-50){\line(0,1){10}}
\multiput(85,-35)(0,5){3}{\circle*{1}}
\multiput(65,-35)(0,5){3}{\circle*{1}}
\multiput(85,-50)(0,10){2}{\circle*{2}}

\end{picture}
\end{center}

\medskip
\noindent {\bf Remark}
If we consider a branch $k[[x_1,\ldots,x_n]]/P_j$ of an algebroid
curve, the sequence of the local blowups of such  a branch does not appear in
the blowup tree of the curve, but the multiplicity sequence of the branch
can be read off in the multiplicity tree of the curve, moving upwards
along the $j$-th branch of the tree and taking the $j$-th component of each
node.
In fact  
the blowup $R_1$ of $R$ has 
exactly one minimal prime, say $Q_j$ contracting to $P_j$
(it is the intersection of the corresponding minimal prime
$k[[t_1]]\times \dots \times k[[t_{j-1}]]\times (0)\times k[[t_{j+1}]]
\times \dots \times k[[t_d]]$ of $\bar R$
with $R_1$) and $R/P_j \subset R_1/Q_j \subseteq k[[t_j]]=
\overline{R/P_j}$. Now the blow up of $R/P_j$ is exactly 
$R_1/Q_j$. 
Moreover, since $n_j \cap R_1$ is the unique maximal 
prime of $R_1$ containing $Q_j$,
we have that $R_1/Q_j=(R_1)_{n_j \cap R_1}/Q_j$.

\medskip
It follows by the previous remark 
that the multiplicity tree of an algebroid curve 
is determined by the
multiplicity sequences of its branches and by the position of the
branching nodes of the tree. In the multiplicity tree of the example above
there is a branching node at level 1.
 
We define two algebroid curves to be {\it equivalent} 
if they have the same number of branches
and the branches can be ordered in a way that gives the same multiplicity
tree. 

Of course this definition extends the equivalence between algebroid
branches recalled in previous section. Let's see how it is also coherent
with the classical equivalence definition for plane curves.

Consider two branches of an algebroid curve $R=k[[x_1,\ldots,x_n]]/ P_1
\cap\cdots\cap P_d$, corresponding to the indices $j$ and
$h$, $1 \leq j,h \leq d$, $j \neq h$. Then $R_{jh} =
 k[[x_1,\ldots,x_n]]/ P_j \cap  P_h = R/\bar P_j \cap \bar P_h$ is an
algebroid curve with two branches and its multiplicity tree is given by
the multiplicity tree of $R$, cancelling all the branches except the
$j$-th and the $h$-th and all the coordinates of the vectors, except the
$j$ and $h$-th. 

In case of a plane algebroid curve,
the {\it intersection number} of two branches is defined as the
intersection number of the ring $R_{jh}$, i.e. $\mu_{jh} = \mu(R/\bar P_j
,R/ \bar P_h)=l_{R_{jh}}(R_{jh}/\bar P_j + \bar P_h)$. By a theorem of Max
Noether (cf. for example \cite [Section 8.4, Theorem 13]{BK}), 
 we have (with our notation):
 \[\mu_{jh}=\sum_{0 \leq i \leq s}  e_i^j e_i^h\]
where $s$ is the level of the branching node of the tree.
Thus the intersection number $\mu_{jh}$ says for plane curves at which
level the two branches split. 

Thus two plane
algebroid curves  $R=k[[x_1,x_2]]/ P_1
\cap\cdots\cap P_d$ and  $T=k[[x_1,x_2]]/ Q_1
\cap\cdots\cap Q_d$ are equivalent (in our sense) if and only if, 
after renumbering the indices, for each
$j$, $1 \leq j \leq d$, $R/\bar P_j$ and $T/\bar Q_j$ have the same
multiplicity sequence and for each $j$ and $h$, $1 \leq j,h \leq d$, $j
\neq h$, $\mu(R/\bar P_j,R/\bar P_h)=\mu(T/\bar Q_j,T/\bar Q_h)$. This is
exactly the classical definition of Zariski of equivalence between plane
algebroid curves (cf. \cite{zar} and \cite{wa}).

\smallskip
\noindent{\bf Remark} Notice that for non plane algebroid curves the
numbers $\mu_{jh}$ does not indicate the level where the $j$-th
and $h$-th branches split. In the previous example $l_R(R/P_1 + P_2)=3$,
but \[\sum_{0 \leq i \leq 1} e_i^1 e_i^2 = 5.\]

\medskip
Let $R=k[[x_1,\ldots ,x_n]]/P_1\cap \cdots \cap P_d$ be an algebroid
curve with $d$ branches and let $S=v(R)$ be its value semigroup. If 
${\boldsymbol\alpha} \in \mathbb N^d$, set $R({\boldsymbol\alpha})=\{r \in
R; v(r)
\geq  {\boldsymbol\alpha}\}$. As in the one-branch case, $R$ is an {\it Arf
ring} if
$x^{-1}R(\boldsymbol\alpha)$, where $v(x)=\boldsymbol\alpha$, is a ring for each 
$\boldsymbol\alpha \in
v(R)$ (equivalently $R$ is Arf if each integrally closed ideal is
stable, cf. for the equivalence \cite [Lemmas 3.18 and 3.22]{BDF}). There
is also in this more general case a smallest Arf overring $R'$ of $R$,
called the Arf closure of $R$ (cf. \cite [Proposition-Definition
3.1]{L}), that has the same multiplicity tree as $R$ (cf. \cite
[Proposition 5.3]{BDF}). Thus $R$ is equivalent to its Arf closure $R'$.

A good semigroup $S \subseteq \mathbb N^d$ is called an {\it Arf semigroup}
if $S(\boldsymbol\alpha)-\boldsymbol\alpha$ is a semigroup, for each ${\bf
\boldsymbol\alpha}
\in S$, where 
$S({\boldsymbol\alpha})= \{\boldsymbol\beta \in S; \boldsymbol\beta \geq {\boldsymbol\alpha}\}$ (equivalently
if each semigroup ideal $S(\boldsymbol\alpha)$ is stable, cf. \cite [p.233]{BDF}). To
each (local) Arf semigroup $S \subseteq \mathbb N^d$ is associated a
multiplicity tree (cf. \cite [Section 5, p.247]{BDF}) and, if $S=v(R')$,
the multiplicity tree of the ring $R$ is the same as the multiplicity
tree of the semigroup $S$ (cf. \cite [Proposition 5.10]{BDF}). Also in
this case the multiplicity tree of a local Arf semigroup characterizes the
semigroup completely and a tree {\bf T} $=\{{\bf e}_{(i)}^j\}$ of vectors
of
$\mathbb N^d$ is the multiplicity tree of a local Arf semigroup if and
only if it satisfies the following conditions:

a) There exists $n\in{\mathbb N}$ such that, for $m\ge n$,
${\bf e}_{(m)}^j=(0,\ldots,0,1,0,\ldots,0)$ (the nonzero coordinate
in the $j$-th position) for any $j=1,\ldots,d$.\\
\indent b) The $h$-th coordinate of ${\bf e}_{(i)}^j$ is $0$ if and only
if ${\bf e}_{(i)}^j$ is not in the 
$h$-th branch of the tree (the $h$-th branch of the tree is
the unique maximal path containing the $h$-th unit vector) and ${\bf
e}_{(i)}^{j_1}\equiv{\bf e}_{(i)}^{j_2}$ (i.e. the two vectors give the
same node in the tree) if and only if the $j_1$-th and $j_2$-th branches
are glued in a node at level $i$.\\
\indent c) ${\bf e}_{(i)}^j=\sum_{ {\bf e}\in {\bf T'}\setminus{\bf
e}_{(i)}^j}{\bf e}$,   for some finite subtree
${\bf T'}$ of {\bf T}, rooted in ${\bf e}_{(i)}^j$\\
(cf. \cite [Theorem 5.11]{BDF}).

We will call such a tree a {\it multiplicity tree} of $\mathbb N^d$.

\medskip
As in the one branch case, we can get the Arf semigroup from the multiplicity
tree taking ${\bf 0}$ and sums of vectors lying on subtrees rooted in the
root of our multiplicity tree. Conversely we can get the multiplicity tree
from the Arf semigroup in the following way: let $S_j$ be the projection
of the Arf semigroup $S$ on the $j$-th coordinate. $S_j$ is a numerical Arf
semigroup (cf. \cite[Proposition 3.30]{BDF}). Denote by $\{e_i^j\}_{i \geq
0}$ its multiplicity sequence. The multiplicity tree is determined by these
multiplicity sequences and by the fact that the $j$-th and $h$-th branches
are glued together as long as the projection of
$S(\boldsymbol\alpha)-\boldsymbol\alpha$, with
${\bf
\alpha}=(e_0^j,e_0^h)+ \dots +(e_k^j,e_k^h)$, on the \lq \lq $jh$-plane"
is local. 
\bigskip
\bigskip
\bigskip
\bigskip
\bigskip
\bigskip
\bigskip
\bigskip
\bigskip
\bigskip\bigskip\bigskip\bigskip\bigskip\bigskip\bigskip\bigskip

\noindent {\bf Example}

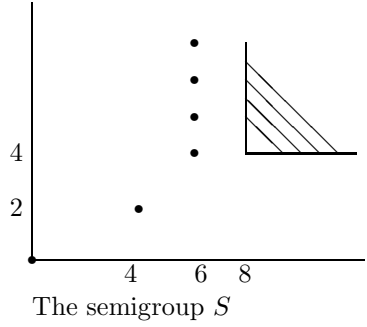
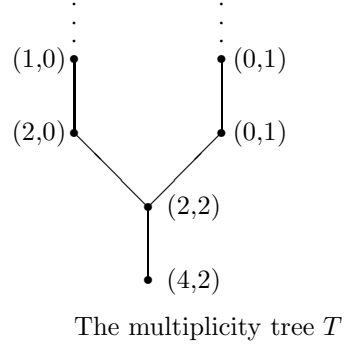
\setlength{\unitlength}{0.49mm}
\begin{center}
\begin{picture}(90,90)(51,-20)

\put(0,0){\line(1,0){90}}
\put(0,0){\line(0,1){70}}

\put(29,14){\circle*{2}}
\put(0,0){\circle*{2}}
\multiput(44,29)(0,10){4}{\circle*{2}}

\put(58,29){\line(1,0){30}}
\put(58,29){\line(0,1){30}}
\put(63,29){\line(-1,1){5}}
\put(68,29){\line(-1,1){10}}
\put(73,29){\line(-1,1){15}}
\put(78,29){\line(-1,1){20}}
\put(83,29){\line(-1,1){25}} 
\put(25,-6){4}
\put(44,-6){6}
\put(56,-6){8}
\put(-6,12){2}
\put(-6,27){4}
\put(0,-15){The semigroup $S$}
\put(150,0){\circle*{2}}
\put(150,0){\line(0,1){20}}
\put(155,-2){(4,2)}
\put(150,20){\circle*{2}}
\put(155,18){(2,2)}
\put(150,20){\line(-1,1){20}}
\put(150,20){\line(1,1){20}}
\multiput(130,60)(40,0){2}{\circle*{2}}
\multiput(130,40)(40,0){2}{\circle*{2}}
\put(130,40){\line(0,1){20}}
\put(113,38){(2,0)}
\put(113,58){(1,0)}
\put(170,40){\line(0,1){20}}
\multiput(173,38)(0,20){2}{(0,1)}
\multiput(170,65)(0,5){3}{\circle*{1}}
\multiput(130,65)(0,5){3}{\circle*{1}}
\put(130,-15){The multiplicity tree $T$}
\end{picture}
\end{center}

The semigroup $S$ contains for example the vectors $(4,2)$,
$(6,5)=(4,2)+(2,2)+(0,1)$, $(9,5)=(4,2)+(2,2)+(0,1)+(2,0)+(1,0)$, obtained
summing vectors along subtrees of $T$.

By projecting $S$ on the two coodinates, we get the multiplicity sequences
$4,2,2,1,\dots$ and $2,2,1, \dots$ respectively. $S(4,2)-(4,2)$ is local, but
$S(6,4)-(6,4)$ is not local, so the tree splits here.

\medskip

We identify two Arf semigroups $S$ and $\widetilde S$ of $\mathbb N^d$ 
(and the corresponding multiplicity trees) if there exists a permutation
$\sigma$ on
$\{1,
\ldots,d\}$ such that $\boldsymbol
\alpha=(e_0^1+\cdots+e_{k_1}^1,\ldots,e_0^d+\cdots+e_{k_d}^d) \in S$ if
and only if $\sigma(\boldsymbol
\alpha)=(e_0^{\sigma(1)}+ \cdots +e_{k_1}^{\sigma(1)},\ldots,e_0^{\sigma(d)}+
\cdots+e_{k_d}^{\sigma(d)})\in \widetilde S$.
 
In conclusion we have:

\begin{thm}\label{bij} The equivalence
classes of algebroid curves with $d$ branches are in one-to-one
correspondence with the Arf  semigroups of $\mathbb
N^d$.
\end{thm}\label{1-1}

\noindent {\bf Proof.} Two algebroid curves with $d$ branches  $R$ and $T$
are equivalent if and only if their Arf closures $R'$ and $T'$ are
equivalent (cf. \cite [Proposition 5.3]{BDF}). This means that the
branches of $T'$ can be ordered in a way such that the multiplicity trees of
$R'$ and $T'$ are the same, i.e. $R'$ and $T'$ have the same value
semigroup, which is an Arf semigroup of $\mathbb N^d$.  On the other hand
each Arf semigroup
$S \subset \mathbb N^d$ is the value semigroup of an algebroid curve with
$d$ branches (cf. \cite[Corollary 5.8]{BDF}).  

\medskip

As in the one branch case, given an equivalence class of algebroid curves 
(equivalentely given an Arf semigroup of $\mathbb N^d$ or a multiplicity tree of
vectors of $\mathbb N^d$), not always there is a plane curve in that class:  a
characterization of multiplicity trees of $\mathbb N^d$  representing
a class that contains a plane curve is easily given, using \cite[Theorem 5]{Bay}. 
As a matter of fact, as we recalled in Section 2, for a plane branch $R$,
knowing the semigroup
$S=v(R)$ is equivalent to knowing the equivalence class of $R$, i.e. the
multiplicity sequence of $R$. Moreover there is a characterization for
numerical semigroups and for multiplicity sequences admissible for plane
branches (cf. e.g. \cite[Proposition 4.8 and Theorem 3.2]{piane}) and, by
\cite[Theorem 5]{Bay}, given two semigroups
$S $ and
$\widetilde S$ admissible for plane branches, the possible
intersection numbers of two plane branches 
with value semigroups $S$ and $\widetilde S$ respectively are known.

Assume we have an ordered set of $d$ multiplicity sequences admissible for
plane branches and that the intersection number between every pair of
consecutive branches is admissible according to \cite[Theorem 5]{Bay}. We
claim that the multiplicity tree obtained in this way is the multiplicity
tree of a plane curve. In fact it follows from the proof of \cite[Theorem
5]{Bay} that, given a plane branch $R$ with $v(R)=S$ and a semigroup
$\widetilde S$ admissible for plane branches, we can find a plane branch
$\widetilde R$ such that $v(\widetilde R)=\widetilde S$ and the
intersection number between
$R$ and
$\widetilde R$ is
$\mu$, where $\mu$ is any of the intersection numbers determined in
\cite[Theorem 5]{Bay}.

\medskip\noindent{\bf Example}
Suppose we have a curve 
with two branches with semigroups $\langle 4,6,13\rangle$
and $\langle 2,3\rangle$, respectively. 
By \cite [Theorem 5]{Bay}
the possible intersection numbers
are $8,12,13$. Since the branches have multiplicity
sequences $4,2,2,1,\ldots$ and $2,1,\ldots$, respectively (cf. e.g. 
\cite[proof of Theorem 2.5]{piane}), it follows
that they split on level 0, 2, or 3.
So we have the following three multiplicity trees. Below any tree is given
a plane curve with that multiplicity tree.
 \setlength{\unitlength}{0.49mm}
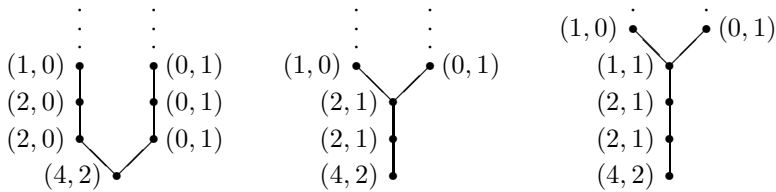
\begin{center}
\begin{picture}(150,70)(76,-90)

\put(55,-72){$(4,2)$}
\put(65,-60){\line(0,1){10}}
\put(85,-60){\line(0,1){10}}
\put(75,-70){\circle*{2}}
\put(75,-70){\line(1,1){10}}
\put(75,-70){\line(-1,1){10}}
\put(45,-62){$(2,0)$}
\put(88,-62){$(0,1)$}
\multiput(45,-52)(0,10){1}{$(2,0)$}
\multiput(45,-42)(0,10){1}{$(1,0)$}
\multiput(88,-52)(0,10){2}{$(0,1)$}
\multiput(65,-60)(20,0){2}{\circle*{2}}
\put(30,-90){$k[[(t^4,u^3),(t^6+t^7,u^2)]]$}

\multiput(65,-50)(0,10){2}{\circle*{2}}
\put(65,-50){\line(0,1){10}}
\put(85,-50){\line(0,1){10}}
\multiput(85,-35)(0,5){3}{\circle*{1}}
\multiput(65,-35)(0,5){3}{\circle*{1}}
\multiput(85,-50)(0,10){2}{\circle*{2}}
 
\put(130,-72){$(4,2)$}
\put(150,-70){\line(0,1){20}}
\multiput(150,-70)(0,10){3}{\circle*{2}}
\put(150,-50){\line(1,1){10}}
\put(150,-50){\line(-1,1){10}}
\multiput(130,-62)(0,10){2}{$(2,1)$}
\multiput(120,-42)(0,10){1}{$(1,0)$}
\multiput(163,-42)(0,10){1}{$(0,1)$}
\multiput(140,-40)(20,0){2}{\circle*{2}}
\put(112,-90){$k[[(t^4,2u^2),(t^6+t^7,u^3)]]$}
 
\put(205,-72){$(4,2)$}
\put(225,-70){\line(0,1){30}}
\multiput(225,-70)(0,10){4}{\circle*{2}}
\put(225,-40){\line(1,1){10}}
\put(225,-40){\line(-1,1){10}}
\multiput(205,-62)(0,10){2}{$(2,1)$}
\multiput(205,-42)(0,10){1}{$(1,1)$}
\multiput(238,-32)(0,10){1}{$(0,1)$}
\put(195,-32){$(1,0)$}
\multiput(215,-25)(20,0){2}{\circle*{1}}
\multiput(215,-30)(20,0){2}{\circle*{2}}
\put(199,-90){$k[[(t^4,u^2),(t^6+t^7,u^3)]]$}

\multiput(160,-35)(0,5){3}{\circle*{1}}
\multiput(140,-35)(0,5){3}{\circle*{1}}
 
\end{picture}
\end{center}


\section{Arf characters of an algebroid curve}

 Our aim is to find, in analogy with the
one-branch case, a finite minimal Arf system of generators for the local
Arf semigroup $v(R')$, i.e. a set of  characters of the algebroid curve
$R$.

We recall some terminology. We call a {\it multiplicity tree} with $d$
branches a tree satisfying conditions a), b), c) given in the previous
section and a {\it
multiplicity branch} a multiplicity tree with one branch, i.e. a
multiplicity sequence of an Arf numerical semigroup.

Suppose that $E$ is a collection of $d$ multiplicity branches $\{e_i^1
\}_{i \geq 0}$, $\ldots$, $\{e_i^d\}_{i\geq 0}$. Denote by
$\tau(E)$ the set of all multiplicity trees having those $d$ branches and
by $\sigma(E)$ the set of the corresponding Arf semigroups.

Consider on $\tau(E)$ the order relation defined setting, for $T_1, T_2
\in \tau(E)$, $T_1 \leq T_2$ if and only if $S(T_1) \subseteq S(T_2)$,
where $S(T_1)$ (resp. $S(T_2)$) is the Arf semigroup defined by $T_1$
(resp. $T_2$).

Before giving the next lemma we need some definitions and notation.

By {\it pinching} a (multiplicity) tree once, we mean modifying a tree
identifying two nodes immediately over a branching node, where the label
of the new node is the coordinatewise sum of the labels of the identified
nodes.

\medskip
\noindent{\bf Example} The following picture describes pinching of a tree
on level 3, level 2, and level 3 consecutively.

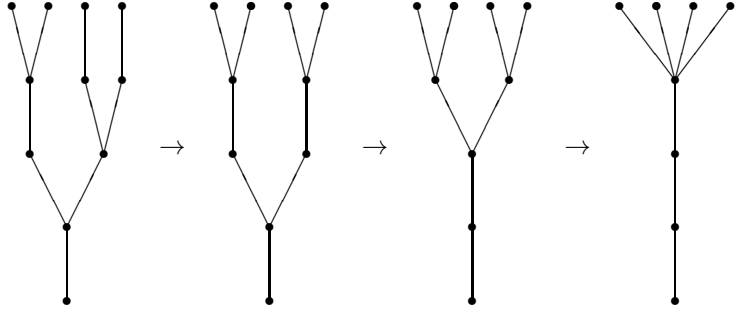
\setlength{\unitlength}{0.49mm}
\begin{center}
\begin{picture}(90,115)(51,-20)

\put(20,0){\circle*{2}}
\put(20,0){\line(0,1){20}}
\put(20,20){\circle*{2}}
\put(20,20){\line(1,2){10}}
\put(20,20){\line(-1,2){10}}
\put(30,40){\circle*{2}}
\put(10,40){\circle*{2}}
\put(30,40){\line(1,4){5}}
\put(30,40){\line(-1,4){5}}
\put(10,40){\line(0,1){20}}
\put(10,60){\line(-1,4){5}}
\put(10,60){\line(1,4){5}}
\put(10,60){\circle*{2}}
\put(15,80){\circle*{2}}
\put(5,80){\circle*{2}}
\put(35,60){\circle*{2}}
\put(25,60){\circle*{2}}
\put(35,60){\line(0,1){20}}
\put(25,60){\line(0,1){20}}
\put(35,80){\circle*{2}}
\put(25,80){\circle*{2}}

\put(45,40){$\rightarrow$}
\put(75,0){\circle*{2}}
\put(75,0){\line(0,1){20}}
\put(75,20){\circle*{2}}
\put(75,20){\line(1,2){10}}
\put(75,20){\line(-1,2){10}}
\put(85,40){\circle*{2}}
\put(65,40){\circle*{2}}
\put(85,40){\line(0,1){20}}
\put(65,40){\line(0,1){20}}
\put(65,60){\circle*{2}}
\put(65,60){\line(1,4){5}}
\put(65,60){\line(-1,4){5}}
\put(85,60){\circle*{2}}
\put(60,80){\circle*{2}}
\put(70,80){\circle*{2}}
\put(70,80){\circle*{2}}
\put(85,60){\line(1,4){5}}
\put(85,60){\line(-1,4){5}}
\put(90,80){\circle*{2}}
\put(80,80){\circle*{2}}

\put(100,40){$\rightarrow$}
\put(130,0){\circle*{2}}
\put(130,0){\line(0,1){20}}
\put(130,20){\circle*{2}}
\put(130,20){\line(0,1){20}}
\put(130,40){\circle*{2}}
\put(130,40){\line(1,2){10}}
\put(130,40){\line(-1,2){10}}
\put(120,60){\circle*{2}}
\put(120,60){\line(1,4){5}}
\put(120,60){\line(-1,4){5}}
\put(140,60){\circle*{2}}
\put(115,80){\circle*{2}}
\put(125,80){\circle*{2}}
\put(125,80){\circle*{2}}
\put(140,60){\line(1,4){5}}
\put(140,60){\line(-1,4){5}}
\put(145,80){\circle*{2}}
\put(135,80){\circle*{2}}

\put(155,40){$\rightarrow$}
\put(185,0){\circle*{2}}
\put(185,0){\line(0,1){20}}
\put(185,20){\circle*{2}}
\put(185,20){\line(0,1){20}}
\put(185,40){\circle*{2}}
\put(185,40){\line(0,1){20}}
\put(185,60){\circle*{2}}
\put(185,60){\line(3,4){15}}
\put(185,60){\line(-3,4){15}}
\put(170,80){\circle*{2}}
\put(180,80){\circle*{2}}
\put(180,80){\circle*{2}}
\put(185,60){\line(1,4){5}}
\put(185,60){\line(-1,4){5}}
\put(200,80){\circle*{2}}
\put(190,80){\circle*{2}}

\end{picture}
\end{center}

Let $T$ be a tree in $\tau(E)$ and let $N$ be a level where all branches
of $T$ are distinct (there is such a level by definition of multiplicity
tree). Then $T$ is determined by $T^{(N)}=(n_1, \ldots, n_{d-1})$, where
$2n_j$ is the distance between the nodes on $j$-th and $j+1$-th branch at
level $N$ (the distance between two nodes is the shortest walk between
them in the tree).

\begin{lem}\label{equiv} Let $T_1, T_2 \in \tau(E)$. Then the following are
equivalent:

(1) $T_1 \leq T_2$;

(2) $T_1$ may be obtained from $T_2$ by a finite number of pinchings;

(3) $T_1^{(N)}$ is coordinatewise less than or equal to $T_2^{(N)}$, where
$N$ is such that $T_1$ and $T_2$ have both $d$ distinct branches at level
$N$. 
\end{lem}

\noindent {\bf Proof.} That (1) $\Leftrightarrow$ (2) follows from the
correspondence between Arf semigroups and multiplicity trees. That (2)
$\Leftrightarrow$ (3) follows from the observation that pinching a tree $T$
once means subtracting 1 from one of the coordinates in $T^{(N)}$.
\medskip

\begin{lem}\label{intersection} If $S_1, S_2 \in \sigma(E)$, then $S_1
\cap S_2 \in \sigma(E)$.
\end{lem}

\noindent {\bf Proof.} We argue on the corresponding multiplicity trees
$T_1, T_2 \in \tau(E)$, showing that, with respect to the order relation
defined above, $\inf(T_1,T_2)$ exists in $\tau(E)$. Let $N \in \mathbb N$
such that $T_1, T_2$ both have $d$ distinct branches at level $N$. If
$T_1^{(N)}=(n_1^1, \ldots, n_{d-1}^1)$ and
$T_2^{(N)}=(n_1^2,
\ldots, n_{d-1}^2)$, we want to show that the tree $T$ described by
$T^{(N)}=(\inf(n_1^1,n_1^2), \ldots, \inf(n_{d-1}^1,n_{d-1}^2))$ is a
multiplicity tree of $\tau(E)$. Suppose that $T \notin \tau(E)$. So there
is in $T$ a node ${\bf e}_{(i)}^j \neq \sum_{{\bf e} \in T' \setminus {\bf
e}_{(i)}^j}{\bf e}$, for each finite subtree $T'$ of $T$, rooted  in ${\bf
e}_{(i)}^j$. This means that
there are two nonzero coordinates $e_i^{j_1}$ and $e_i^{j_2}$ of ${\bf
e}_{(i)}^j$ such that 
$$e_i^{j_1}=e_{i+1}^{j_1}+e_{i+2}^{j_1}+ \cdots+ e_{i+k_1}^{j_1}$$
$$e_i^{j_2}=e_{i+1}^{j_2}+e_{i+2}^{j_2}+ \cdots +e_{i+k_2}^{j_2}$$
with $k_1 < k_2$ and that the $j_1$-th and $j_2$-branches split at a level
 strictly greater than $i+k_1$. It follows that in $T_1$ or in $T_2$ the
$j_1$-th and $j_2$-th branches are glued at level $i$ and split at a level
 strictly greater than $i+k_1$. This means that $T_1$ or $T_2 \notin
\tau(E)$, a contradiction.

\medskip

Suppose we have an Arf semigroup $S$ in $\mathbb N^d$ or, equivalently, a
multiplicity tree with $d$ branches. We will describe a method to find
a finite set $V$ of vectors in $S$ which determines $S$.

For $1 \leq j \leq d $, let $\{e_i^j\}_{i \geq 0}$ be  the multiplicity
branches of the tree, let $S_j = \{e_0^j, e_0^j + e_1^j, \ldots\} $ be the
Arf numerical semigroup of the $j$-th branch and let
$\{c_1^j,\cdots,c_{n_j}^j\}$ be the Arf system of generators for $S_j$.
To construct $V$, consider, for $1 \leq j \leq d$ and for $1 \leq i \leq
n_j$, the vectors of
$S$ of the following form:
 $$(\alpha_1, \ldots, \alpha_{j-1}, c_i^j, \alpha_{j+1}, \ldots,
\alpha_d)$$
 where, fixed $c_i^j$, the other coordinates are minimal (this choice is
possible by property (2) of good semigroups (cf. Section 4)). Each of these
vectors corresponds in $T_S$, the multiplicity tree associated to $S$, to a
subtree with only one branch, i.e.  to a node.
For any couple of consecutive branches (the $j$-th and the $j+1$-th
branch) of $T_S$ there is a node in the tree where they split. Let's call
it the
$(j,j+1)$- branching node of the tree. For each $j$, $1 \leq j \leq d-1$,
check if in the set $V$ there is a vector corresponding to a node of the
tree over the
$(j,j+1)$-branching node. If not, add to the set $V$ such a vector.

\begin{thm}\label{V} The smallest Arf semigroup containing $V$ is $S$.
\end{thm}

\noindent {\bf Proof.} By our choice, all the numbers $c_1^j, \cdots,
c_{n_j}^j$ appear as $j$-th coordinate of some of the vectors in $V$. So
the numerical semigroups
$S_j={\rm Arf}(c_1^j, \cdots, c_{n_j}^j)$ (and the collection $E$ of their
multiplicity sequences) are determined. It follows that any Arf semigroup
containing $V$ is in $\sigma(E)$ and there exists at least one such (the
Arf semigroup corresponding to a tree with a unique branching node at
level 0). 

Notice that the presence in an Arf semigroup of $\sigma(E)$ of a vector
${\bf v}=(\alpha_1,\ldots,\alpha_d)$, with $\alpha_j=e_0^j+ \cdots +
e_{k_j}^j$ for $1 \leq j \leq d$, and with $k_{j_0}>k_{{j_0}+1}$ for some
${j_0}$, forces the ${j_0}$-th and ${j_0}+1$-th branches of the
corresponding tree to split at most at level
$k_{{j_0}+1}$.  
By the choice of the vectors in $V$, any Arf semigroup containing $V$
corresponds to a multiplicity tree of $\tau(E)$ with $d $ distinct
branches at level $N$, where $N$ is the greatest $k_j$ which appears in
the expressions ${\bf v}=(\alpha_1,\ldots,\alpha_d)$, with $\alpha_j=e_0^j+
\cdots +
e_{k_j}^j$ for $1 \leq j \leq d$, ${\bf v}$ in $V$. It follows that
the number of Arf semigroups  of $\sigma(E)$ containing $V$ is finite. 
By Lemma \ref{intersection}, the intersection $H$ of all these is still an Arf
semigroup of $\sigma(E)$ and so it is the smallest Arf semigroup
containing $V$. Let $T_H$ be the tree corresponding to $H$ and set
$T_H^{(N)}=(n_1,\dots,n_d)$ (cf. notation before Lemma \ref{equiv}). Let
${j_0} \in \{1,\dots,d-1\}$ and let ${\bf v}=(\alpha_1,\ldots,\alpha_d)$,
with $\alpha_j=e_0^j+ \cdots +
e_{k_j}^j$ for $1 \leq j \leq d$ be a vector of $V$ corrsponding to a node
above the $({j_0},{j_0}+1)$-branching node of $T_S$. We have $k_{j_0} \neq
k_{{j_0}+1}$. Supposing $k_{j_0} > k_{{j_0}+1}$, i.e. supposing the node of
${\bf v}$ on the
${j_0}$-th branch, we have, by the minimality of $H$, $n_{j_0} =
N-k_{{j_0}+1}$. On the other hand, since ${\bf v}$ correspond to a node
above the $({j_0},{j_0}+1)$ branching node of $T_S$, the ${j_0}$-th and
${j_0}+1$-th branches of
$T_S$ are forced to split exactly at level $k_{{j_0}+1}$. This means that
setting
$T_S^{(N)}=(m_1,\dots,m_{d-1})$, we have also $m_{j_0} = N-k_{{j_0}+1}$. It
follows that $T_H=T_S$ and so $H=S$.

\medskip

\noindent {\bf Remark} The set $V$ is not uniquely determined, as is
clear by its construction. 
We can e.g. get a smaller set by eliminating vectors in $V$ 
with the following rule:
if ${\bf v_1}$, ${\bf v_2}$ and ${\bf v_3}=\min({\bf v_1},{\bf v_2})\neq
{\bf v_i}$ (for $i=1,2$) are in $V$,
then eliminate ${\bf v_3}$. 

\medskip

But even for minimal subsets of $V$ determining $S$ 
the cardinality is not unique. 
A very simple counterexample is the following. 
Let $$S=\{(0,0,0), (1,1,1)\}\cup ((2,2,2)+\mathbb N^3)$$  
Then
$V_1=\{(1,1,1),(2,3,2)\}$ and $V_2=\{(1,1,1), (3,2,2), (2,2,3)\}$ both are
minimal subsets determining $S$.

In any case
it seems proper to call a minimal set $V$ determining $S$ 
a set of {\it characters} of $R$, where $v(R')=S$.

From the discussion above the following result follows:

\begin{thm}\label{subrings} Let $R$ be an algebroid curve which is an Arf
ring, with
$v(R)=S$, and let
 $\{\bf v_1, \ldots ,
\bf v_N\}$ be a set of characters of $R$. Then any subring $U=k[[f_1,
\dots, f_N]]$ of
$R$ with $v(f_i)=\bf v_i$, $i=1,\dots, N$, has Arf closure  $R$ and so $U$
is an algebroid curve equivalent to $R$.
\end{thm}

In some cases, for a particular Arf ring $R$, a proper subset of 
$\{f_1, \dots, f_N\}$ is enough to generate a ring with Arf closure $R$.
In that case the curve $k[[f_1, \dots, f_N]]$
can be projected into a space of lower dimension without changing 
its multiplicity tree, i.e. without changing its equivalence class.

\medskip

We conclude with two examples.

\medskip\noindent
{\bf Example} We consider the following multiplicity tree

 \setlength{\unitlength}{0.49mm}
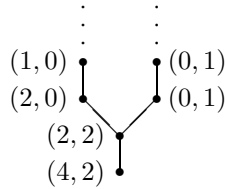
\begin{center}
\begin{picture}(110,55)(26,-75)

\put(55,-72){$(4,2)$}
\put(55,-62){$(2,2)$}
\multiput(45,-52)(0,10){1}{$(2,0)$}
\put(45,-42){$(1,0)$}
\multiput(88,-52)(0,10){2}{$(0,1)$}
\multiput(75,-70)(0,10){2}{\circle*{2}}
\put(75,-70){\line(0,1){10}}
\put(75,-60){\line(-1,1){10}}
\put(75,-60){\line(1,1){10}}
\multiput(65,-50)(0,10){2}{\circle*{2}}
\put(65,-50){\line(0,1){10}}
\put(85,-50){\line(0,1){10}}
\multiput(85,-35)(0,5){3}{\circle*{1}}
\multiput(65,-35)(0,5){3}{\circle*{1}}
\multiput(85,-50)(0,10){2}{\circle*{2}}

\end{picture}
\end{center}

\noindent
which corresponds to the Arf semigroup $$\{(0,0),(4,2),(6,n);n\ge4\}
\cup((8,4)+{\mathbb N}^2).$$
The Arf rings over $k$ with that semigroup (that are all equivalent by
Theorem \ref {bij}) are the rings 
$$k+(t^4+\alpha_5t^5+\alpha_7t^7,
u^2+\beta_3u^3)k+(t^6+\gamma_7t^7,0)k+t^8k[[t]] \times u^4k[[u]],$$
where $\alpha_5,\alpha_7,\beta_3,\gamma_7\in k$.
 The characters for the first branch are $4,6,9$, and for the
second $2,5$. The vectors containing these characters are
$(4,2),(6,4)=(4,2)+(2,2),(9,4)=(4,2)+ (2,2)+(2,0)+(1,0)$ and
$(6,5)=(4,2)+(2,2)+(0,1)$. We can delete $(6,4)=
\min(\{(9,4),(6,5)\}$ and get a minimal set of Arf generators $\{(4,2),(9,4),(6,5)\}$.
For the choice of parameters $\alpha_5=\alpha_7=\beta_3=\gamma_7=0$, we get the ring
$$R=k+(t^4,u^2)k+(t^6,0)k+t^8k[[t]] \times u^4k[[u]]$$
 Any subring of $R$ generated by elements of values
$(4,2),(9,4)$, and $(6,5)$, such as
$$U=k[[(t^4,u^2),(t^9,u^4),(t^6,u^5)]]=k[[x,y,z]]/(x^3-z^2,x^2z-y^2)
\cap(y-x^2,xy^2-z^2)$$ has Arf closure $R$, i.e. is equivalent to $R$. With
another choice of parameters,
$\alpha_5=\alpha_7=\beta_3=0,\gamma_7=1$, we get the ring 

$$\widetilde R=k+(t^4,u^2)k +(t^6+t^7,0)k+t^8k[[t]] \times u^4k[[u]]$$ 

Already the ring $$\widetilde U=k
[[(t^4,u^2),(t^6+t^7,u^5)]]=k[[x,y]]/(y^4-x^7+x^6-4x^5y-2x^3y^2)\cap(x^5-y^2)$$
has Arf closure $\widetilde R$,
so any ring 
$\widetilde U_r =k[[(t^4,u^2),(t^6+t^7,u^5),r]]$ with
$v(r)=(9,4)$ is equivalent to $\widetilde R$ and has a projection to a
plane curve equivalent to
$\widetilde R$.

\medskip\noindent
{\bf Example} We consider the following multiplicity tree

 \setlength{\unitlength}{0.49mm}
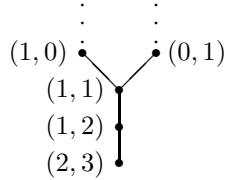
\begin{center}
\begin{picture}(150,55)(76,-75)

\put(130,-72){$(2,3)$}
\put(150,-70){\line(0,1){20}}
\multiput(150,-70)(0,10){3}{\circle*{2}}
\put(150,-50){\line(1,1){10}}
\put(150,-50){\line(-1,1){10}}
\multiput(130,-62)(0,10){1}{$(1,2)$}
\put(130,-52){$(1,1)$}
\multiput(120,-42)(0,10){1}{$(1,0)$}
\multiput(163,-42)(0,10){1}{$(0,1)$}
\multiput(140,-40)(20,0){2}{\circle*{2}}
\multiput(160,-37)(0,5){3}{\circle*{1}}
\multiput(140,-37)(0,5){3}{\circle*{1}}
 
\end{picture}
\end{center}

\noindent
which corresponds to the Arf semigroup
$$\{(0,0),(2,3),(3,5)\}\cup((4,6)+{\mathbb N}^2)$$ The Arf rings over
$k$ with that semigroup are the rings 
$$k+(t^2,
u^3+\beta_4u^4+\beta_5u^5)k+(t^3,\gamma u^5)k+t^4k[[t]] \times u^6k[[u]],
\gamma\ne0,$$
where $\beta_4,\beta_5,\gamma\in k$.
The characters for the first branch are $2,3$, and for the second $3,5$.
The vectors containing these characters are $(2,3)$ and $(3,5)=(2,3)+(1,2)$, so we have to
add one more generator corresponding to a node above the branching point, say
$(4,7)=(2,3)+(1,2)+(1,1)+(0,1)$.
If we choose the parameters $\beta_4=\beta_5=0,\gamma=1$, we can take the ring
$$k[[(t^2,u^3),(t^3,u^5),(t^4,u^7)]]=k[[x,y,z]]/(z-x^2,xz-y^2)
\cap(yz-x^4,xz-y^2,x^3y-z^2)$$ as an example of a ring with
$$k+(t^2,u^3)k+(t^3,u^5)k+t^4k[[t]] \times u^6k[[u]]$$ as Arf closure.
The choice of a node above the branching point is arbitrary. If we choose
$(6,6)= (2,3)+(1,2)+(1,1)+(1,0)+(1,0)$ we can, with the same choice of
parameters, take
$$k[[(t^2,u^3),(t^3,u^5),(t^6,u^6)]]=k[[x,y,z]]/(x^3-y^2,z-y^2)
\cap(z-x^3,y^3-xz^2)$$ as an example of a ring with
$$k+(t^2,u^3)k+(t^3,u^5)k+t^4k[[t]] \times u^6k[[u]]$$ as Arf closure.

\end{document}